\documentclass[12pt]{article}
\usepackage{amsmath}
\newtheorem{theorem}{Theorem}

\title{The distribution of longest run lengths in integer compositions}
\author{Herbert S. Wilf}
\begin{document}
\maketitle
\begin{abstract}
We find the generating function for $C(n,k,r)$, the number of
compositions of $n$ into $k$ positive parts all of whose runs
(contiguous blocks of constant parts) have lengths less than $r$,
using recent generalizations of the method of Guibas and Odlyzko for
finding the number of words that avoid a given list of subwords.
\end{abstract}
\section{Introduction}
A composition of an integer $n$ is a representation $n=a_1+a_2+\dots+a_k$ in which the parts $a_i$ are positive integers, and where the order of the parts is important. Thus $9=1+1+1+4+2$ is one of the compositions of $n=9$, and $9=4+1+1+2+1$ is another.

A \textit{run} in a composition is a maximal string of consecutive
identical parts. The composition \[28=3+5+5+5+3+3+4\] has run
lengths of 1,3,2,1, for example. In this note we find (the
generating function of) $C(n,k,r)$, the number of compositions of
$n$ into $k$ parts, whose runs all have lengths $<r$ (see Theorem
\ref{th:rungen} below), by using recent generalizations of the
Guibas-Odlyzko theory of counting words that avoid a given list of
subwords.

In their 1981 paper \cite{go}, Guibas and Odlyzko gave an elegant solution to the following counting problem. Given an alphabet ${\cal A}$, and a list ${\cal L}$ of words over that alphabet, the list being \textit{reduced} in the sense that no word on the list is a subword of any other. How many words of length $n$ do not contain any of the words in ${\cal L}$ as a subword? Other solutions of this problem have been given by the cluster method of Goulden and Jackson \cite{gj}, and by Zeilberger's \cite{zz} method of counting words that avoid ``mistakes.''

The results of \cite{go} have recently been extended by A.N. Myers \cite{anm} to the situation wherein the letters of the alphabet are assigned weights, the weight of a word is the sum of the weights of its letters, and one is to find the number of words of weight $n$ that avoid the members of the list ${\cal L}$. This allows us to solve problems involving compositions of integers as well as problems that do not involve compositions.

Finally, Myers's results have been complemented by Heubach and Kitaev \cite{hk} to provide the number of words of \textit{length $k$ and weight $n$} that avoid the members of the list ${\cal L}$, though their theorems are restricted to the alphabet $\{1,2,\dots,n\}$ and therefore apply almost exclusively to integer compositions.

The above theorems present the generating function for the desired numbers of words as the first component of the solution vector of a system of linear, simultaneous equations, or, by using Cramer's rule, as a ratio of two determinants.

The main point of this note is the following. The easy case in such
word problems is the case in which every pair of distinct words on
the forbidden list ${\cal L}$ has correlation 0, in a sense to be
explained below, or equivalently, for every pair $x,y$ of distinct
words on that list, no suffix of $x$ is also a prefix of $y$. In
that situation, the matrix of coefficients of the system of linear
equations that expresses the answer to the question has a very
simple form. It consists of a nonzero first row and first column and
main diagonal, all other entries being 0's.

For a matrix of that form it is easy to write out the solution of
the governing system of linear equations simply and explicitly. We
will do that below and then find the generating function for
$C(n,k,j)$, the number of compositions of $n$ into $k$ parts the
lengths of whose runs is at most $j$.
\section{The main theorem}
Let $X$ and $Y$ be two words over a given alphabet. We define \textit{the correlation $c_{XY}$ of $X$ on $Y$}, as follows.
\begin{itemize}
\item Write the word $X$ above the word $Y$, aligned so that the rightmost letter of $X$ is above the rightmost letter of $Y$.
\item Fix some integer $j\ge 0$. Shift $Y$ $j$ places to the left, so the rightmost letter of $Y$ is now under the $(j+1)$st letter of $X$, counting from the right.
    \item Examine the subword of $X$ that now overlaps with $Y$. This is the maximal prefix of $X$ that has letters of the shifted $Y$ below it.
    \item If that subword of $X$ is identical with the subword of $Y$ that lies below it, take $c_j=1$, else take $c_j=0$.
        \item Having done this for all $j$, the correlation of $X$ on $Y$ is the binary vector $c_0c_1c_2\dots$.
\end{itemize}
For example, if $X=110$ and $Y=1011$ then $c_{XY}=011$ and $c_{YX}=0010$, in which we have written the bits of the $c$'s in the order $c_0c_1\dots c_{m-1}$.

Let each letter $u$ of the alphabet be assigned a weight $w(u)$, and let the weight of a word be the sum of the weights of its letters. Finally, if $X$ is an $m$-letter word $X=a_0a_1\dots a_{m-1}$, define the \textit{correlation polynomial} $c_{XY}(x,q)$ of $X$ on $Y$ to be
\begin{equation}
\label{eq:corr}
c_{XY}(x,q)=
c_0+c_1x^{w(a_{m-1})}q+c_2x^{w(a_{m-2}a_{m-1})}q^2
+\dots+c_{m-1}x^{w(a_1a_2\dots a_{m-1})}q^{m-1}.
\end{equation}
The main result of \cite{hk}, which extends the main result of \cite{anm}, which in turn extends the main result of \cite{go}, is the following.
\begin{theorem}[Heubach, Kitaev]
\label{th:hk}
Let ${\cal L}=\{S_1,\dots,S_k\}$ be a list of integer compositions, such that no composition  on the list is contained in any other. Let $F(x,q)=\sum_{\sigma}x^{w({\sigma})}q^{\ell(\sigma)}$, the sum being extended over all compositions of all integers that avoid every word on the list ${\cal L}$, where  $\ell(\sigma)$ is the length of the word (number of parts of) $\sigma$ and $w({\sigma})$ is the sum of the parts of $\sigma$. Then $F(x,q)$ is the component $x_1$ of the solution vector of the following system of linear equations:
\begin{equation}
\label{eq:simul}
\left(\begin{array}{ccccccc}
1-x(1+q)&1-x&\dots&1-x\\
x^{w(S_1)}q^{\ell(S_1)}&-c_{11}(x,q)&\dots&-c_{1k}(x,q)\\
\vdots&\vdots&\ddots&\vdots\\
x^{w(S_k)}q^{\ell(S_k)}&-c_{k1}(x,q)&\dots&-c_{kk}(x,q)
\end{array}\right)\left(
\begin{array}{ccc}
x_1\\x_2\\\vdots\\x_{k+1}
\end{array}
\right)=\left(\begin{array}{ccc}
1-x\\0\\\vdots\\0
\end{array}\right)
\end{equation}
\end{theorem}
\section{The easy case}
We now specialize to the case where $c_{ij}(x,q)=0$ for all $i\neq j$, $1\le i,j\le k$, $k$ being the length of the forbidden word list ${\cal L}$. The coefficient matrix entries in the equations (\ref{eq:simul}) then all vanish except for those in the first row, the first column, and the main diagonal. For any such matrix, $B$, say, the first entry of the solution vector of the equations $B\mathbf{x}=(1-x,0,\dots,0)^T$ is easily verified to be
\[x_1=\frac{1-x}{b_{11}-b_{12}\frac{b_{21}}{b_{22}}-\dots
-b_{1,k+1}\frac{b_{k+1,1}}{b_{k+1,k+1}}}.\]
If we apply this result to the equations (\ref{eq:simul}) we obtain
\begin{theorem}
\label{th:easy} Let ${\cal L}=\{S_1,\dots,S_k\}$ be a list of
integer compositions, such that no word on the list is contained in
any other. Suppose further that all correlation polynomials
$c_{ij}(x,q)=0$, for $i\neq j$. Let
$F(x,q)=\sum_{\sigma}x^{w({\sigma})}q^{\ell(\sigma)}$, the sum being
extended over all compositions $\sigma$ that avoid every word on the
list ${\cal L}$, where $\ell(\sigma)$ is the length of the word
$\sigma$. Then we have the explicit formula
\begin{equation}
\label{eq:fexp}
F(x,q)=\frac{1}{1-\frac{qx}{1-x}+\sum_{j=1}^k\frac{x^{w(S_j)}q^{\ell(S_j)}}{c_{j,j}(x,q)}}.
\end{equation}
\end{theorem}
\section{Carlitz compositions and beyond}
\subsection{Carlitz compositions}
We apply the results of the previous section to finding the
distribution function of the lengths of the longest runs of integer
compositions. Again, a run in a composition is a maximal string of
identical parts. The composition 28=3+5+5+5+3+3+4 has run lengths of
1,3,2,1, for example.

A Carlitz composition is one all of whose runs have length 1. That is, a Carlitz composition
is one in which no two consecutive parts are equal. These compositions have been extensively
 studied in recent years, both exactly and asymptotically \cite{car, kp, lp}. The machinery
 of ``the easy case'' above counts Carlitz compositions of $n$, as follows.

The list ${\cal L}$ of forbidden subwords is ${\cal L}=\{11,22,33,44,\dots\}$. A Carlitz
 composition is evidently one that avoids this list, and also evidently, this list belongs
  to the easy case, i.e., the off-diagonal correlation polynomials all vanish. Thus we can
   use Theorem \ref{th:easy}.

The word $S_j$ is $jj$, and its weight is $w(S_j)=2j$. The correlation polynomials $c_{S_iS_j}$
 vanish for all $i\neq j$, while for $i=j$ we have by (\ref{eq:corr}),
\[c_{S_jS_j}(x,y)=1+x^jq.\]
If $C(n,k)$ is the number of Carlitz compositions of $n$ into $k$ parts, we now have
 from equation (\ref{eq:fexp}),
\begin{eqnarray*}
\label{eq:cargen}
\sum_{n,k}C(n,k)x^nq^k&=&\frac{1}{1-\frac{xq}{1-x})+q^2\sum_{j\ge 1}\frac{x^{2j}}{1+qx^j}}\\
&=&1+qx+qx^2+(q+2q^2)x^3+(q+2q^2+q^3)x^4+(q+4q^2+2q^3)x^5+\dots
\end{eqnarray*}
This generating function has previously been found, in somewhat different form, by Knopfmacher
 and Prodinger \cite{kp}.
\subsection{Beyond}
Now we find the distribution function of the maximum run length in compositions of $n$ that have $k$ parts.

Let $C(n,k,r)$ denote the number of compositions of $n$ into $k$ parts that have no run of length $\ge r$. Note that $C(n,k,2)$ counts Carlitz compositions of $n$ with $k$ parts. To find $C(n,k,r)$ we use the list ${\cal L}=\{1^r,2^r,3^r,\dots\}$ of forbidden words, where, e.g., $1^r$ is a string of $r$ 1's. Then again the list ${\cal L}$ qualifies for ``the easy case,'' since the correlations all vanish off of the diagonal. while on the diagonal,
\[c_{S_jS_j}(x,y)=1+x^jq+x^{2j}q^2+\dots+x^{j(r-1)}q^{r-1}=\frac{1-q^rx^{rj}}{1-qx^j}.\]
We now have from equation (\ref{eq:fexp}),
\begin{theorem}
\label{th:rungen}
The number $C(n,k,r)$ of compositions of $n$ into $k$ parts that have no run of length $\ge r$ has the generating function
\begin{equation}
\label{eq:rungen}
\sum_{n,k}C(n,k,r)x^nq^k=\frac{1}{1-\frac{xq}{1-x}+q^r\sum_{j\ge 1}\frac{x^{rj}(1-qx^j)}{1-q^rx^{rj}}}.
\end{equation}
\end{theorem}
When $r=3$ we have
\[\sum_{n,k}C(n,k,3)x^nq^k=1+qx+(q+q^2)x^2+(q+2q^2)x^3+(q+3q^2+3q^3)x^4+\dots,\]
and for $r=4$,
\[\sum_{n,k}C(n,k,4)x^nq^k=1+qx+(q+q^2)x^2+(q+2q^2+q^3)x^3+(q+3q^2+3q^3)x^4+\dots.\]
The average length of the longest run in a composition of $n$ has been found to be $\sim \log_2{n}$, by Grabner et al \cite{gkp}, using the method of i.i.d. geometric random variables.


\newpage

\begin{thebibliography}{aaa}
\bibitem{car} L. Carlitz, Restricted compositions, The Fibonacci Quart., \textbf{14} (1976), 254–264.
\bibitem{gj} I. P. Goulden and D. M. Jackson, An inversion theorem for cluster decompositions of sequences with distinguished subsequences,
J. London Math. Soc. \textbf{20}  (1979), no. 3, 567--576.
\bibitem{gkp} P. Grabner, A. Knopfmacher and H. Prodinger, Combinatorics of geometrically distributed random variables: Run statistics, Theoretical Computer Science \textbf{297} (2003), 261--270.
\bibitem{go} L.J. Guibas and A.M. Odlyzko, String overlaps, pattern matching, and nontransitive games, J. Combinatorial Theory, Ser. A, \textbf{30} (1981), 183--208.
\bibitem{hk} Silvia Heubach and Sergey Kitaev, Avoiding substrings in compositions,      \texttt{arXiv:math/0903.5135} [math.CO]
\bibitem{hm} S. Heubach and T. Mansour, Enumeration of 3-letter patterns in compositions,    \texttt{arXiv:math/0603285v1} [math.CO].
\bibitem{kp} Arnold Knopfmacher and Helmut Prodinger, On Carlitz compositions,  European J. Combin.  \textbf{19}  (1998),  no. 5, 579--589.
    \bibitem{lp} Guy Louchard and Helmut Prodinger,  Probabilistic analysis of Carlitz compositions,  Discrete Math. Theor. Comput. Sci.  \textbf{5}  (2002),  no. 1, 71--95.
\bibitem{anm} Amy N. Myers, Forbidden substrings on weighted alphabets, Australasian J. Math., to appear.
    \bibitem{zz} Doron Zeilberger, Enumeration of words by their number of mistakes.  Discrete Math.  \textbf{34}  (1981), no. 1, 89--91.
\end{thebibliography}
\end{document}